\def\IS{{\Bbb S}} 
 
\def\IZ{{\Bbb Z}} 
 
\def\IC{\Bbb C} 
\def\ID{{\Bbb D}}
\def\IA{{\Bbb A}}
\def\oC{\hat{\IC}}

\documentclass[12pt]{article} 

\usepackage[psamsfonts]{amssymb} 
\usepackage{amsfonts,amsmath} 
\usepackage{graphicx}

\newtheorem{theorem}{Theorem}
\newtheorem{corollary}{Corollary}

\usepackage{amsmath, amssymb, amscd, amsxtra, txfonts}

\title{Quasicircles as equipotential lines,  homotopy classes and geodesics.} 
\author{Gaven   J. Martin \thanks{Research supported by the Marsden Fund.  \newline 
Mathematics Subject Classification Primary 30C62, 37F30, 30C75}} 
\date{} 
\begin{document}

\maketitle

\begin{abstract}  We give an application of our earlier results concerning the quasiconformal extension  of a germ of a conformal map to establish that in two dimensions the equipotential level lines of a capacitor are quasicircles whose distortion depends only on the capacity and the level.  As an application we find that given disjoint, nonseparating and nontrivial continua $E$ and $F$ in $\oC=\IC\cup\{\infty\}$,  the closed hyperbolic geodesic generating the fundamental group $\pi_1\big(\oC\setminus (E\cup F) \big) \cong \IZ$ is a $K$-quasicircle separating $E$ and $F$ with explicit distortion bound depending only on the capacity of $\oC\setminus (E\cup F)$.  This result is then extended to obtain distortion bounds on a quasicircle representing a given homotopy class of a simple closed curve in a planar domain.  Finally we are able to use these results to show that a simple closed hyperbolic geodesic in a planar domain  is a quasicircle  with a distortion bound depending explicitly, and only, on its length.\end{abstract}
    
\section{Introduction.}

A {\em capacitor},  sometimes called a {\em condenser}, is simply the complement of a pair of disjoint continua $E,F\subset \IC$.  We assume $E$ and $F$ that $E$ and $F$ do not separate in $\IC$ (and so are {\em cellular continua}).  The conformal capacity of  $\Omega = \oC\setminus E\cup F$ is defined as 
\begin{equation}\label{capdef}
Cap(\Omega) =    \inf_{u}\; \; \frac{1}{4\pi }\;  \iint_\Omega |\nabla u|^2 \; dz  
\end{equation}
where the infimum is taken over all admissible functions $u$,  this means that for $z\in E$, $u(z)=+1$ and for $z\in F$,  $u(z)=-1$. It is well known that there is a unique extremal function $u$ which attains the infimum for the variation problem at (\ref{capdef}).  The physical interpretation is to consider charging $E$ and $F$ with opposite total charge (although this does not matter) which equidistributes on $E$ and $F$.  The extreme function $u$ for (\ref{capdef}) is harmonic and the maximum principle shows $-1< u < 1$ and the level lines (or circles) of $u$ are the equipotential lines of the electric field.  As a well known example the capacity of the annulus $A=\IA(r,R)=\{z:r< |z| < R\}$ is  $Cap(A) =\log R/r$ and the extremal function is  
$u(z)=2 \frac{\log |z|/r}{\log R/r} -1$. \\
{ 
\scalebox{0.35}{\includegraphics[viewport= -180 420 250 790]{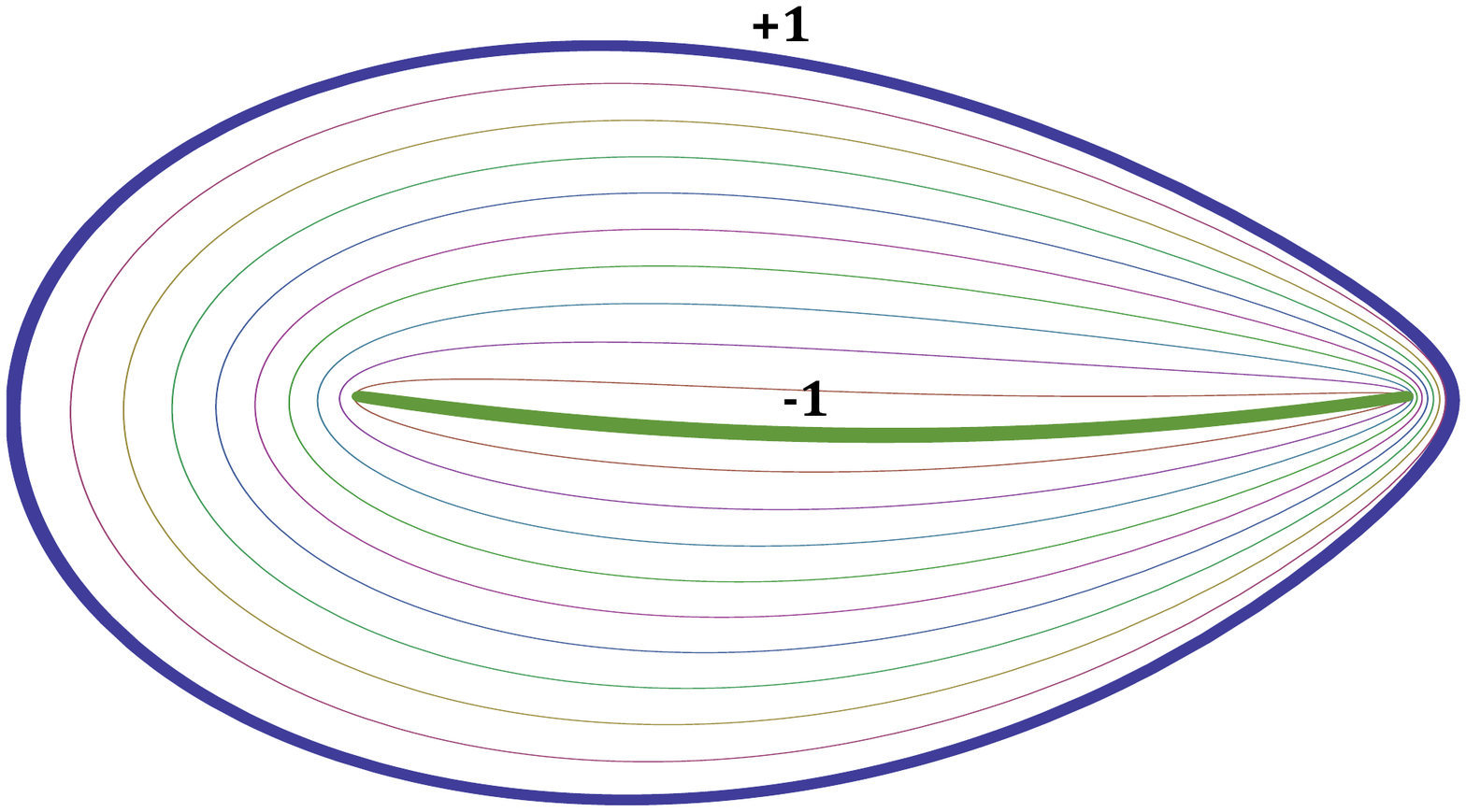}}
}  \\
\noindent{\bf Fig. 1}  Typical equipotential lines of a domain.

\medskip

If $E$ and $F$ are finite points (and when suitably interpreted) the line $u=0$ is the perpendicular bisector and all other level lines are circles. But it is not always the case that some equipotential line is a circle.  The main result of this article is to establish that this is always the case up to a global deformation of bounded geometry - that is a quasiconformal mapping.  We now define these mappings. 

\medskip

 A homeomorphism $f:\Omega\to \IC$  in the Sobolev class $f\in W^{1,2}_{loc}(\Omega,\IC)$ of functions with locally square  integrable first derivatives is called {\em quasiconformal} if there is a $1\leq K <\infty$ for which $f$ satisfies the {\em distortion inequality} 
\begin{equation}\label{de}
|Df(z)|^2 \leq K \; J(z,f), \hskip20pt\mbox{almost everywhere in $\Omega$}
\end{equation}
Here  $Df(z)$ is the Jacobian matrix and $J(z,f)$ its determinant.  If such a $K$ exists,  we will say that $f$ is  {\em $K$-quasiconformal.} The modern theory of quasiconformal mappings can be found in \cite{AIM}.  Such mappings have the basic property of ``bounded distortion'' as they take infinitesimal circles to infinitesimal ellipses - conformal mappings have no distortion ($K=1$).  The distortion inequality (\ref{de}) actually implies the improved regularity $f\in W^{1,2K/(K+1)}_{loc}(\Omega)$,  \cite[Astala's theorem]{AIM}.

\medskip

A {\em $K$-quasicircle}  is the image of the unit circle under a quasiconformal mapping of $\IC$.  These Jordan curves have very many nice geometric properties -  for instance the Ahlfors bounded turning condition.  A Jordan curve $\alpha$ is a quasicircle if and only if there is a constant $C$ with the following property. Each pair of points $z,w\in \alpha$ separates $\alpha$ into two components, $\alpha^+$ and $\alpha^{-}$ and 
\[ \max \Big\{ {\rm diam}(\alpha^+),{\rm diam}(\alpha^-) \Big\} \leq C |z-w|  \]
There are many other equally important criteria and these are catalogued in the book of Gehring and Hag, \cite{GH}.

\section{Results.}

\begin{theorem}\label{main}  Let $\Omega=\oC\setminus E\cup F$ be a capacitor and $u$ the extremal function for (\ref{capdef}).  Then the level curve $\{z:u(z) = a\}$ is a $K_a$-quasicircle and
\begin{equation}
K_a \leq \frac{1}{1-|a|} \;  \left(1+\frac{8\beta_0}{Cap(\Omega)}\right)  \hskip15pt \beta_0 = 2.4984 \ldots
\end{equation}
\end{theorem}

There is a clear consequence here.  Given a simple closed curve $\alpha$ in $\Omega$ we can define 
\[
Cap_\Omega(\alpha)=\sup  \Big\{  Cap(D) : \mbox{ $D \subset \Omega$ is doubly connected and $\alpha$ generates $\pi_i(D)$ }\Big\}
\]
It is not difficult to see that $Cap_\Omega(\alpha)=\infty$ if and only if $\alpha$ is either homotopically trivial,  or homotopic to an isolated boundary point. In each case there is a round circle in the homotopy class of $\alpha$ in $\Omega$.  If $Cap_\Omega(\alpha)<\infty$,  we can apply Theorem \ref{main} to the $a=0$ level line for any doubly connected domain $D\subset \Omega$ as above and use the elementary compactness properties of quasiconformal mappings to obtain the following theorem.

\begin{theorem}\label{thm2} Let $\alpha$ be a simple closed curve in a planar domain $\Omega$.  Then there is a $K$-quasicircle homotopic to $\alpha$ with
\begin{equation}
K  \leq    1+\frac{8\beta_0}{Cap_\Omega(\alpha)} 
\end{equation}
\end{theorem}
\medskip
There is an obvious candidate for the quasicircle given to us by Theorem \ref{thm2} and that is the hyperbolic geodesic in the homotopy class of $\alpha$,  should $\Omega$ admit a hyperbolic metric.  It is generically true that for planar domains that there there is a simple  closed  hyperbolic geodesic in the homotopy class of a nontrivial simple closed curve, the exception occurring when there are isolated points in the boundary.  We are able to establish the following theorem.

\begin{theorem}\label{new3}  Let $\Omega$ be a planar domain with a hyperbolic metric and $\alpha$ a simple closed geodesic of length $\ell$.  Then $\alpha$ is a $K$-quasicircle and 
\[ K \leq 1+ \frac{4 \beta_0\, \ell }{\pi  \arccos \Big(\tanh \Big(\frac{\ell}{2}\Big)\Big) } \]
\end{theorem}
This estimate is nontrivial for all lengths $\ell\geq 0$ and has nice behaviour as $\ell\to 0$,  but is exponential in $\ell$ (roughly $e^{\ell/2}$) for large $\ell$.  This behaviour seems incorrect and perhaps an artefact of the method of proof here which is based on collaring theorems for torsion free Fuchsian groups.  For doubly connected domains $\Omega$,  Theorem \ref{main} and Theorem \ref{thm2},  together with the length formula at (\ref{length}) below,  give a better bound.

\begin{theorem}  Let $\Omega$ be a doubly connected planar domain with a hyperbolic metric and $\alpha$ the simple closed geodesic of length $\ell$.  Then $\alpha$ is a $K$-quasicircle and 
\begin{equation}
K  \leq   1+\frac{4\,\beta_0 \, \ell }{\pi^2} \end{equation}
\end{theorem}

\section{Example.}
Before proving Theorem \ref{main} it is tempting to discuss a conjectured extremal mapping.  In view of the extremal mapping for quasicircles identified in \cite{MarProc} which approximates the situation here when one of the continua degenerates we expect that the extremal domain is a twisted Teichm\"uller domain given by deleting vertical and horizontal segments.  Thus set 
\[ \Omega_0 = \IC \setminus \Big( \Big\{|\Im m(z)|\geq \frac{1}{2} \Big(m^{-2}- m^2), \Re e(z)=0\Big\} \cup \Big\{|\Re e(z)|\leq 1, \Im m(z)=0 \Big\} \Big) \]
The uniformising mapping $\IA(r,1/r) \mapsto \Omega$ is a hyperbolic isometry effected by the conformal transformations
\begin{eqnarray*}
\IA(r,1/r) &  \stackrel{z\mapsto rz}{\rightarrow} & \IA(1,r^2) \stackrel{z\mapsto \frac{1}{2}(z+1/z) }{\rightarrow} {\cal E}_{r} \setminus [-1,1] \to \ID \setminus [-m^2,m^2]  \\
& \stackrel{z\mapsto \frac{1}{z}}{\rightarrow} & \IC \setminus \big( |x|\geq 1/m^2 \cup \ID\big)\stackrel{z\mapsto \frac{1}{2}(z-1/z) }{\rightarrow} \Omega
\end{eqnarray*}
Here $ {\cal E}_r$ is the ellipse for which the sum of the lengths of the major axis and minor axis is $r^2$.  The map ${\cal E}_{r}\to \ID$ is a Riemann map,  given by
\[ z\mapsto J_{SN}\Big(\frac{2}{\pi} \, {\cal K}(m^2)\arcsin(\frac{z}{m}) \big| m^2\Big)\]
where $J_{SN}$ is  the Jacobi elliptic sine function and 
\[ \frac{\pi}{4} \;\frac{{\cal K}'(m)}{{\cal K}(m)} = \log 1/r^2 \]
Here ${\cal K}$ and ${\cal K}'$ are the elliptic integrals of the first and second kind.  This is illustrated below when $r=\sqrt{2}$.\\
{ 
\scalebox{0.4}{\includegraphics[viewport= -70 120 450 500]{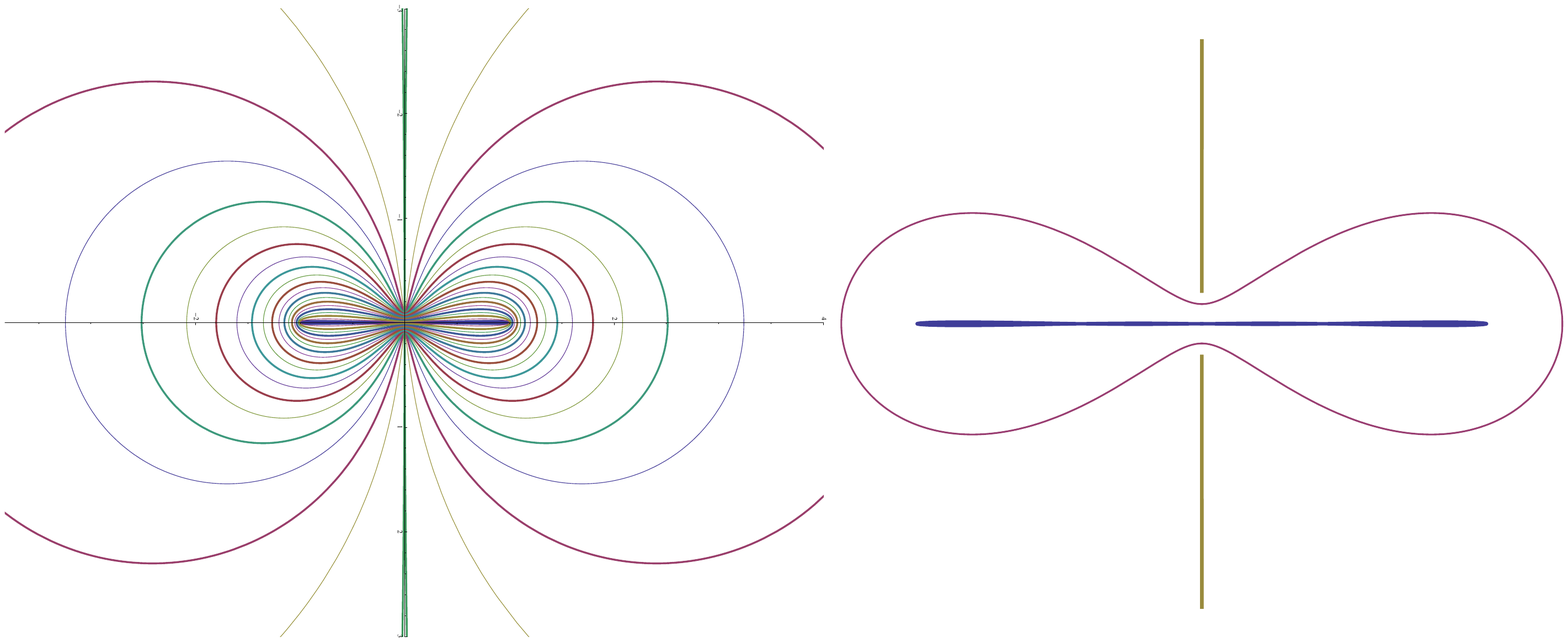}}
}  \\
\noindent{\bf Fig. 2} A conjectured extremal configuration with $\Omega$ of conformal modulus $\log 2$ and the associated level curve $h=0$ (scaled). This curve is the hyperbolic geodesic generating the fundamental group of $\Omega_0$.

\section{Proof of Theorem \ref{main}.} 
The domain $\Omega = \oC\setminus E\cup F$ is doubly connected,  has a single conformal invariant and is therefore conformally equivalent to the annulus $\IA(r,1/r)=\{z:r<|z|<1/r\}$ where $r$ is chosen so that the  conformal capacity of $\Omega$ is 
$Cap(\Omega) = \log r^2$.  
The extremal function for $\IA(r,1/r)$ is 
\begin{equation}
u(z)=  \frac{\log\big( |z|^2/r^2\big)}{\log\big( 1/r^2\big)}  -1
\end{equation}
and so the extremal function for $\Omega$ is $v = u(\varphi)$,  where $\varphi:\Omega\to \IA(r,1/r)$ is the (unique up to rotation) conformal mapping.

To complete the  proof   we will subsequently establish the following two theorems.

\begin{theorem}\label{thm5} Let $v$ be the extreme function for the capacitor $\Omega$.  Then the equipotential line $\{v=0\}$ is a $K_0$-quasicircle and 
\[ K_0 \leq  1+\frac{8\beta_0}{Cap(\Omega)}  \hskip15pt \beta_0 = 2.4984 \ldots  \]
\end{theorem}
\begin{theorem}\label{thm6} Let $v$ be the extreme function for the capacitor $\Omega=\oC\setminus (E\cup F)$.  Then for $-1<a\leq b<1$,  there is a $K_{ab}$--quasiconformal mapping $f:\oC\to\oC$ with $f|E=idenity$ and $f|F=identity$ and which takes the  equipotential line $\{v=a\}$ to the equipotential line $\{v=b\}$ and has  distortion 
\[ K_{ab} = \max \left\{ \frac{b+1}{a+1} , \frac{1-b}{1-a} \right\} .  \]
\end{theorem}

\medskip

Now if we put $b=0$ we obtain $f:\oC\to\oC$ with $f(\{v=0\})= \{v=a\}$ and $f$ is quasiconformal with distortion $K=\max \left\{ \frac{1}{a+1} , \frac{1}{1-a} \right\}=\frac{1}{1-|a|}$.   Since $\{v=0\}$ is a quasi circle with bound given by Theorem \ref{thm2} the result will follow.

 \section{Proof of Theorem \ref{thm5}}  This proof is based on some arguments in earlier work with Klen (see   \cite{KM})  concerning the quasiconformal extension of a holomorphic germ and it significantly sharpens an earlier result in \cite{MarProc}.  It is in the spirit of \cite{DM} where David,  Semmes and the author show the euclidean bisecting curve $\{z: d(z,E)=d(z,F) \}$ is a quasicircle.  
 
 The equipotential line $\gamma=\{v=0\}\subset \Omega$ is a smooth Jordan curve (possibly including $\infty$) and bounds two simply connected regions $\Omega_+$,  containing $E$, and $\Omega_-$ containing $F$.  Further $\gamma = \varphi^{-1}(\IS)$.  Let $\psi:\Omega_+\to \ID$ be a Riemann mapping.  Then the mapping $\eta = \psi\circ \varphi^{-1}:\IA(r,1) \to \ID \setminus \psi(E)$ is conformal and has $\eta(\IS)=\IS$ when extended via the Carath\'eodory theorem.  Thus $(\ID\setminus \psi(E),\eta^{-1})$ is a {\em holomorphic germ} for the quasisymmetric mapping $g=\eta|\IS:\IS\to\IS$.  Then according to \cite[Theorem 3.1]{KM} the quasisymmetric mapping $g$ has a $K$-quasiconformal extension $\tilde{\eta}:\ID\to\ID$,  $\tilde{\eta}|\IS=g$, to a self map of the disk with the distortion bound
\begin{equation}\label{kbound} K \leq 1+ \frac{4\beta_0}{\log 1/r},  \hskip15pt \beta_0 = 2.4984 \ldots \end{equation}
 Actually $\beta_0$ is the modulus of the Teichm\"uller ring $\IC\setminus \big( [0,1/\sqrt{2}]\cup [1,\infty) \big)$ and this result is sharp in the sense that the constant $4\beta_0$ cannot be replaced by a number less than $1$.  If more information is known - for instance if the capacity is known to be large - then better estimates are given in \cite{KM} and we will use these later,  see (\ref{later}).  Now $\psi^{-1}\circ \tilde{\eta}:\ID \to \Omega+$ is quasiconformal with the bound at (\ref{kbound}) and 
 \[ \psi^{-1}\circ \tilde{\eta}|\IS =\psi^{-1}\circ  \eta |\IS =  \psi^{-1}\circ  \psi \circ \varphi^{-1} |\IS =  \varphi^{-1}|\IS \] 
 Thus $f_+=\psi^{-1}\circ \tilde{\eta}:\Omega_+ \to \ID$ is quasiconformal with distortion bound at (\ref{kbound}) and agrees with $\varphi^{-1}$ on the boundary.  In the same way we may construct $f_-:\Omega_-\to \oC \setminus \overline{\ID}$ with the same distortion bound and also agreeing with $\varphi^{-1}$ on $\IS$.  It is a moments work to now see that $f:\oC\to\oC$ defined by 
 \begin{equation}
 f(z) = \left\{ \begin{array}{ll} f_+(z) & z\in \ID \\
 \varphi^{-1}(z) & z\in \IS \\f_-(z) & z\in \oC\setminus \overline{\ID}
 \end{array} \right.
 \end{equation}
 is $K$-quasiconformal with distortion bound given at (\ref{kbound}). Finally we note that  
 \[  1+ \frac{4\beta_0}{\log 1/r} =1+ \frac{8\beta_0}{\log 1/r^2} = 1+ \frac{8\beta_0}{Cap(\IA(r,1/r))} = 1+ \frac{8\beta_0}{Cap(\Omega)} \] 
 and this completes the proof. \hfill $\Box$
 
 \medskip
 
When we equip $\Omega$ with the hyperbolic metric,  the level line $\alpha = \{v=0\}$ above is actually the hyperbolic geodesic generating the fundamental group of $\Omega$ as $\varphi$ is a hyperbolic isometry and $\IS$ is the hyperbolic geodesic generating the fundamental group of $\IA(r,1/r)$.  We therefore have the following corollary.
 
 \begin{corollary} \label{cor1} Let  $E$ and $F$ be nontrivial disjoint continua in $\oC=\IC\cup\{\infty\}$.  Endow $\Omega=\oC\setminus (E\cup F) $ with a complete hyperbolic metric.  Then the closed hyperbolic geodesic generating the fundamental group $\pi_1\big(\Omega\big) \cong \IZ$ is a $K$-quasicircle separating $E$ and $F$ with explicit distortion bound  
 \[ K\leq 1+   \frac{8\beta_0}{Cap(\Omega)},  \hskip15pt \beta_0 = 2.4984 \ldots \] 
  \end{corollary}
 Notice that this distortion bound has reasonable behaviour as $Cap(\Omega)\to\infty$ and $Cap(\Omega)\to 0$.  It seems likely that the asymptotics are correct and that the constant $8\beta_0$ cannot be made smaller than $2$.

  \section{An extremal problem}  Here we want to study the following extremal problem for quasiconformal mappings.  Let $0<r<1$ and identify the quasiconformal mapping of smallest distortion for which $f(\IA(r,1/r))= \IA(r,1/r)$ and  for $r < s \leq  t < 1/r$ we require  $f(\IS(0,s))=\IS(0,t)$. We can identify a lower bound by considering various moduli.  The hypotheses imply that $f(\IA(r,s))=\IA(r,t)$ and that $f(\IA(s,1/r))=\IA(t,1/r)$.  Any $K$-quasiconformal mapping can distort capacity by at most a factor $K$ and hence any mapping with these properties must have distortion
\begin{eqnarray}
K & \geq & \frac{Cap(\IA(r,t))}{Cap(\IA(r,s))} = \frac{\log (t/r)}{\log(s/r)} = \frac{\log t-\log r }{\log s - \log r} \label{1} \\
K & \geq & \frac{Cap(\IA(t,1/r))}{Cap(\IA(s,1/r))} = \frac{\log (tr)}{\log (sr) } = \frac{\log t+ \log r }{\log s + \log r} \label{2} 
\end{eqnarray}
  Indeed it is well known that the extremal quasiconformal mappings between round annuli are the power mappings $(r,\theta)\mapsto (r^\alpha,\theta)$,  $\alpha>0$,  after suitable normalisation,  and that these mappings are uniquely extremal up to a rotation.   Thus there is a map which effects both (\ref{1}) and (\ref{2}) and this radial mapping has the same boundary values on the common boundary $\IS(0,t)$ and further,  this mapping is the identity on $\IS(0,r)$ and $\IS(0,1/r)$ so that these bounds are both achieved. This establishes the next theorem.
  
  \begin{theorem} \label{extremal} Let $0<r<1$ and $r<s, t<1/r$.  Then there is a $K$-quasiconformal homeomorphism $f:\IA(r,1/r)\to \IA(r,1/r)$ such that $f|\partial \IA(r,1/r)=identity$,  $f(\IS(0,s))=\IS(0,t)$ and 
  \[ K = \max \left\{ \frac{\log t-\log r }{\log s - \log r} , \frac{\log t+\log r }{\log s + \log r}  \right\} \]
  Such a mapping is unique up to rotation and has minimal distortion among all quasiconformal homeomorphisms irrespective of the boundary values.
  \end{theorem}
  
  \section{Proof of Theorem \ref{thm6}}
  
Recall that the extremal function for $\Omega$ is the composition $v = u(\varphi)$ where $u(z)=  \log( |z|^2/r^2)/\log(1/r^2)  -1$ and  $\varphi:\Omega\to \IA(r,1/r)$ is conformal.  The level curve $v=a$ is therefore mapped by $\varphi$ to the circle where $u=a$ and this is the circle of radius $r^{-a}$.  Next,  by Theorem \ref{extremal} there is a $K$-quasiconformal self mapping $f$ of $\IA(r,1/r)$ which is the identity on the boundary,  takes the circle of radius $r^{-a}$ to the circle of radius $r^{-b}$ and has distortion at most 
\begin{equation}\label{kk}
K = \max \left\{ \frac{\log r^{-b} -\log r }{\log r^{-a} - \log r} , \frac{\log r^{-b}+\log r }{\log r^{-a} + \log r}  \right\} =    \max \left\{ \frac{b+1}{a+1} , \frac{1-b}{1-a} \right\} . 
\end{equation}

Now suppose that the boundaries of $E$ and $F$ are smooth Jordan curves.  Then the map $\varphi$ extends continuously and homeomorphically to the boundary $\partial\Omega = \partial E \cup \partial F \to\partial \IA(r,1/r)$.  Then the mapping $h = \varphi^{-1} \circ f \circ \varphi :\Omega \to \Omega$ (here $f$ is as given by Theorem \ref{extremal}) extends continuously to the boundary as the identity, takes the level curve $\{v=a\}$ to the level curve $\{v=b\}$ and  is $K$-quasiconformal,  with $K$ as at (\ref{kk}).   We can therefore extend this mapping to the extended complex plane by setting $h(z)=z$ for $z\in E\cup F$ and observe that $h$ remains $K$-quasiconformal.   Now an arbitrary non separating continua $E$ can be approximated by a smoothly bounded disk $D_E$,  with $E\subset D_E$,  in the Hausdorff topology.  The family of quasiconformal self maps of $\oC$ that we obtain from the above argument with approximates $D_{E,\epsilon}$ and $D_{F,\epsilon}$ all have a uniform distortion bound,  all are the identity on $E$ and $F$,   and therefore this family does not converge to a constant in $\oC\setminus \{z_0\}$ for any point $z_0\in \oC$.  The compactness properties of  quasiconformal mappings now provides us with a limit $K$-quasiconformal mapping which has all  the desired properties to complete the proof once we further observe that the corresponding conformal maps will converge locally uniformly on $\oC \setminus (E\cup F)$ and so for fixed $a$ and $b$ the limit level curves will be the level curves of the limit.\hfill $\Box$

\section{Hyperbolic Geodesics}   The hyperbolic metric density of the annulus  $\IA_r=\IA(r,1/r)$ can be found in \cite[\S 12.2]{BeardonMinda} as
\begin{equation*}\label{hypden}
d_{\IA_r}(z) = \frac{\pi}{2 \log(1/r)}  \; \frac{1}{|z| \cos\Big(\frac{\pi\log|z|}{2\log(1/r)}  \Big)}.
\end{equation*}
From this formula we compute that the length of the geodesic circle $\{|z|=1\}$ is 
\begin{equation}\label{length} \ell =  \frac{\pi^2}{\log(1/r)}. 
\end{equation}
The length of the geodesic generating the fundamental group of a doubly connected hyperbolic region uniquely determines the conformal equivalence class.  

Next,  we are going to establish a consequence of a result of A. Beardon.  This result is stated in his book in terms of a universal constraint on the geometry of hyperbolic elements in a Fuchsian group and when reinterpreted we find the following theorem.  Hyperbolic metrics are assumed complete and normalised to have constant curvature equal to $-1$. 

\begin{theorem}\label{thmhyp} Let $\Omega$ be a planar domain with hyperbolic distance $\rho_\Omega$.  Let $\alpha$ be a simple closed geodesic in this metric and let $\ell$ be the hyperbolic length of $\alpha$.  Let $r$ and $r_0$ be defined by the equations
\begin{eqnarray*}
r = e^{-\pi^2/\ell},  & &  r_0=e^{-2\pi \arccos \big(\tanh \big(\frac{\ell}{2}\big)\big)/\ell}
\end{eqnarray*}
Then there is a hyperbolic  isometry $h$ defined on the annular subregion $\IA(r_0,1/r_0) \subset \IA(r,1/r)$, with the hyperbolic metric of $\IA(r,1/r)$, into $\Omega$ so that $h(\IS)=\alpha$.
\end{theorem}
\noindent{\bf Proof.}  The uniformisation theorem asserts that up to isometry we can realise $\Omega$ as the orbit space $\ID/\Gamma$ of a torsion free Fuchsian group $\Gamma$ acting as hyperbolic isometries of the hyperbolic plane  $\ID$.  The geodesic $\alpha$ lifts to a hyperbolic line $\tilde{\alpha}$ in $\ID$ which is realised as the axis of a hyperbolic transformation $f$ with translation length $\ell$.  Since $\alpha$ is simple,  the primitive hyperbolic $f$ is simple  and thus no translate of this axis $\tilde{\alpha}$ crosses $\tilde{\alpha}$: that is for each $g\in \Gamma$ either
\[ g(\tilde{\alpha}) = \tilde{\alpha},  \hskip5pt {\rm or} \hskip5pt g(\tilde{\alpha}) \cap \tilde{\alpha} =\emptyset. \]
For each $g\in \Gamma$ the element $f^g=gf g^{-1} $ is hyperbolic,  has translation length  $\ell$ - the same $f$,  and has axis $g(\tilde{\alpha})$.  The group $\langle f,g \rangle$ is nonelementary unless $g\in \langle f\rangle$.  Now  \cite[Corollary 11.6.10]{Beardon} asserts that for  $g\not\in \langle f\rangle$,
\[
\sinh^2\big(\frac{\ell}{2}\big) \cosh(\rho_{hyp}(\tilde{\alpha},g(\tilde{\alpha})) \geq \cosh^2\big(\frac{\ell}{2}\big)+1
\]
In particular we have achieved the bound 
\begin{equation}\label{10} \cosh(\rho_{hyp}(\tilde{\alpha},g(\tilde{\alpha})) \geq 1+\frac{2}{\sinh^2\big(\frac{\ell}{2}\big)}.
\end{equation}
If we choose $\delta_0$ so that $\cosh(2\delta_0) = 1+2/{\sinh^{2}\big(\frac{\ell}{2}\big)}$,  then the bound at (\ref{10}) implies that the neighbourhood ${\cal C}={\cal C}(\alpha,r_0)= \{z\in \ID: \rho_{hyp}(z,\alpha) <\delta_0\}$ has the property that for all $g\in \Gamma$ either $g({\cal C}) = {\cal C}$ or $g({\cal C}) \cap {\cal C} =\emptyset$.  Thus ${\cal C}$ projects down to ${\cal C}/\langle f \rangle$ - an embedded annulus of hyperbolic radius $\delta_0$ about $\alpha$ in $\Omega$.  As noted above,  this annulus can be identified uniquely up to conformal equivalence (or hyperbolic isometry)  as the neighbourhood of hyperbolic radius $r_0$ about $\IS$ in $\IA(r,1/r)$.  We will have proved the theorem once we solve the equation
\begin{equation}
\delta_0 = \frac{\pi}{2 \log(1/r)}  \int_{1}^{1/r_0}    \; \frac{dt}{t \cos\Big(\frac{\pi\log t}{2\log(1/r)}  \Big)}.
\end{equation}
which arises from the fact that the radial lines in $\IA(r,1/r)$ are geodesic and the hyperbolic density formula given at (\ref{hypden}) can be used to calculate distance.
Performing the integral,  this equation reads as
\begin{eqnarray*}
\delta_0 & = &   \log \left[\frac{1+\tan \Big( \frac{\pi}{4} \, \frac{ \log(1/r_0)}{\log(1/r)} \Big)}{1-\tan \Big( \frac{\pi}{4} \, \frac{ \log(1/r_0)}{\log(1/r)} \Big) } \right],\\
  \cosh(2\delta_0) & = & \frac{1}{2} \left( \left( \frac{1+\tan \Big( \frac{\pi}{4} \, \frac{ \log(1/r_0)}{\log(1/r)} \Big)}{1-\tan \Big( \frac{\pi}{4} \, \frac{ \log(1/r_0)}{\log(1/r)} \Big) } \right)^2+\left( \frac{1-\tan \Big( \frac{\pi}{4} \, \frac{ \log(1/r_0)}{\log(1/r)} \Big)}{1+\tan \Big( \frac{\pi}{4} \, \frac{ \log(1/r_0)}{\log(1/r)} \Big) } \right)^2 \right) \\
 & = &  2\sec^2\left( \frac{\pi}{2} \, \frac{ \log(1/r_0)}{\log(1/r)} \right) -1
\end{eqnarray*}
Now,  by our original choice of $\delta_0$ and as $\log(1/r)=\pi^2/\ell$, we are led to the equation
\begin{eqnarray*}
 2\sec^2\left( \frac{\pi}{2} \, \frac{ \log(1/r_0)}{\log(1/r)} \right) -1 & = & 1+\frac{2}{\sinh^{2}\big(\frac{\ell}{2}\big)}, \\ 
  \cos\left( \frac{\pi}{2} \, \frac{ \log(1/r_0)}{\log(1/r)} \right)  & = & \tanh \Big(\frac{\ell}{2}\Big)\\
  \log(1/r_0)  & = &\frac{2\pi}{\ell} \arccos \Big(\tanh \Big(\frac{\ell}{2}\Big)\Big) 
\end{eqnarray*}
  This now completes the proof of the theorem.
 \hfill $\Box$

\medskip

We now come to the following result  which is really the point of the above calculation.

\begin{theorem}  Let $\Omega$ be a planar domain with a hyperbolic metric and $\alpha$ a simple closed geodesic of length $\ell$.  Then $\alpha$ is a $K$-quasicircle and 
\[ K \leq 1+ \frac{4 \beta_0\, \ell }{\pi  \arccos \Big(\tanh \Big(\frac{\ell}{2}\Big)\Big) } \]
\end{theorem}
\noindent{\bf Proof.}  We apply Corollary \ref{cor1} to the annular region about $\alpha$ which is conformally equivalent to $\IA(r_0,1/r_0)$ of capacity $\log 1/r_{0}^{2}$. The result follows for the estimate of $r_0$ given by Theorem \ref{thmhyp}. \hfill $\Box$

\medskip

\noindent{\bf Remark.} Note that the function $ \frac{\pi}{\ell} \arccos \Big(\tanh \Big(\frac{\ell}{2}\Big)\Big)$ is strictly decreasing and convex.  A little calculation then gives the cleaner estimate
\begin{equation}
K \leq 1+ \frac{5\ell }{\pi }\, \sqrt{e^\ell +1}.
\end{equation} 
If $\ell$ is small, so the modulus of the annulus is large,   then \cite[Theorem 3]{KM} actually gives a better bound (which also applies to Corollary \ref{cor1} if the capacity is large).  We need $\log 1/r_0\geq 2 \beta_0$ for this, 
and after a little manipulation and simplification we achieve the estimate 
\begin{equation}\label{later}
K \leq 1+ \frac{3\ell }{2}, \hskip15pt \mbox{whenever $\ell \leq 1$}.
\end{equation}

\bigskip

\noindent G.  Martin -  Massey University,  Auckland, NZ,    g.j.martin@massey.ac.nz

\end{document}